\documentclass[10pt]{article}
\usepackage{amsmath}
\usepackage{amssymb}
\usepackage{latexsym}
\begin{document}
\title{{\bf Hom-Bol algebras}}
\author{{\it Sylvain Attan and A. Nourou Issa}}
\date{}
\maketitle
\begin{abstract}
Hom-Bol algebras are defined as a twisted generalization of (left) Bol algebras. Hom-Bol algebras generalize multiplicative Hom-Lie triple systems in the same way as Bol algebras generalize Lie triple systems. The notion of an $n$th derived (binary) Hom-algebra is extended to the one of an $n$th derived binary-ternary Hom-algebra and it is shown that the category of Hom-Bol algebras is closed under the process of taking $n$th derived Hom-algebras. It is also closed by self-morphisms of binary-ternary Hom-algebras. Every Bol algebra is twisted into a Hom-Bol algebra. Relying on the well-known classification of real two-dimensional Bol algebras, examples of real two-dimensional Hom-Bol algebras are given.
\end{abstract}
\section{Introduction}
A {\it Bol algebra} is a triple $(A, \cdot , [, ,])$, where $A$ is a vector space, $\cdot : A^{\otimes 2}  \rightarrow A$ a bilinear map (the binary operation on $A$) and $[, ,]: A^{\otimes 3}  \rightarrow A$ a trilinear map (the ternary operation on A) such that\\
\par
(B1) $x \cdot y = - y \cdot x$,
\par
(B2) $[x,y,z] = - [y,x,z]$,
\par
(B3) ${\circlearrowleft}_{x,y,z} [x,y,z] = 0$,
\par
(B4) $[x, y, u \cdot v ] = [x,y,u] \cdot v + u \cdot [x,y,v] + [u,v, x \cdot y] - uv \cdot xy$,
\par
(B5) $ [x,y, [u,v,w]] = [[ x,y,u],v,w ] + [u, [x,y,v],w ] + [u,v, [x,y,w] ]$ \\
\\
for all $u,v,w,x,y,z$ in $A$, where ${\circlearrowleft}_{x,y,z}$ denotes the sum over cyclic permutation of $x,y,z$ and juxtaposition will be used (here and in the sequel) in order to reduce the number of braces (so, in (B4), $uv \cdot xy$ means $(u \cdot v) \cdot (x \cdot y)$).
\par
Observe that when $x \cdot y = 0$ in a Bol algebra $(A, \cdot , [, ,])$, then it reduces to a {\it Lie triple system} $(A,[, ,])$ so one could think of a Bol algebra $(A, \cdot , [, ,])$ as a Lie triple system $(A,[, ,])$ with an additional anticommutative binary operation ``$\cdot$'' such that (B4) holds.
\par
The definition of a Bol algebra given above is the one of a {\it left} Bol algebra (see, e.g., \cite{Mikh}, \cite{Perez}, \cite{Sab-Mikh1}) and the reader is advised not to confuse it with the one of the {\it right} Bol algebra $(A, \diamond, (; ,))$ (see, e.g., \cite{Sab}). However, a left Bol algebra $(A, \cdot , [, ,])$ is obtained from a right one $(A, \diamond, ( ; , ))$ if set $x \cdot y = -x \diamond y$, and $[x,y,z] = -(z;x,y)$. In this paper, by a Bol algebra we always mean a left Bol algebra.
\par
Bol algebras are introduced (see \cite{Mikh}, \cite{Sab-Mikh1}, \cite{Sab-Mikh2}) in a study of the differential geometry of smooth Bol loops. Such a study could be seen as a generalization of the differential geometry of Lie groups, where the left-invariant vector fields on a given Lie group constitute a Lie algebra. The well-known correspondence between Lie groups and Lie algebras has been extended to some classes of smooth loops. Such an extension is first performed for Moufang loops by A. I. Mal'tsev \cite{Malt}, and it proved that one has a one-to-one correspondence between {\it Malcev algebras} \cite{Sagle} and Moufang loops (see \cite{Kuz}, \cite{Malt}, \cite{Nagy}). Later, the study of that correspondence has been extended to smooth Bol loops in \cite{Mikh}, \cite{Sab-Mikh1}, \cite{Sab-Mikh3} (see also \cite{Sab} and \cite{Sab-Mikh4} for a survey of the subject; for the infinitesimal theory, including the differential geometry, of local smooth loops, one may refer to \cite{Sab}, \cite{Sab-Mikh2}, \cite{Sab-Mikh4}). The tangent structure to a smooth Bol loop turns out to be a Bol algebra and, locally, there is a correspondence between Bol algebras and (local) smooth Bol loops (\cite{Mikh}, \cite{Sab-Mikh1}). For local smooth loops in general, the correspondence is established in terms of a vector space equipped with a family of multilinear operations, called {\it hyperalgebra} \cite{Sab-Mikh3} (now called {\it Sabinin algebra}). Bol algebras are further studied in \cite{Kuz-Zaid}, \cite{Perez}.
\par
It is shown \cite{Mikh} that Bol algebras, as tangent algebras to local smooth Bol loops, are Akivis algebras with additional conditions. {\it Akivis algebras} are introduced (\cite{Ak1}, \cite{Ak2}) in a study of the differential geometry of 3-webs (see references in \cite{Ak1}, \cite{Ak2}; originally, M. A. Akivis called ``$W$-algebras'' such algebraic structures and the term ``Akivis algebra'' is introduced in \cite{Hofm}). Akivis algebras have also a close connection with the theory of smooth loops \cite{Ak1}.
\par
The aim of this paper is a study of a Hom-type generalization of Bol algebras. Roughly, a Hom-type generalization of a kind of algebra is obtained by a certain twisting of the defining identities by a linear self-map, called the twisting map, in such a way that when the twisting map is the identity map, then one recovers the original kind of algebra. In this scheme, e.g., associative algebras and Leibniz algebras are twisted into Hom-associative algebras and Hom-Leibniz algebras respectively \cite{Makh-Silv1} and, likewise, Hom-type analogues of Novikov algebras, alternative algebras, Jordan algebras or Malcev algebras are defined and discussed in \cite{Makh1}, \cite{Yau3}, \cite{Yau4}. The Hom-type generalization of some classes of ternary algebras are discussed in \cite{Atag}, \cite{Yau3}, \cite{Yau5}. One could say that the theory of Hom-algebras originated in \cite{Hart} in a study of deformations of the Witt and the Virasoro algebras (in fact, some $q$-deformations of the Witt and the Virasoro algebras have a structure of a Hom-Lie algebra \cite{Hart}). Some algebraic abstractions of this study are given in \cite{Makh-Silv1}, \cite{Yau1}, \cite{Yau2}. For more recent results regarding Hom-Lie algebras or Hom-Leibniz algebras, one may refer to \cite{Cas1}, \cite{Cas2}, \cite{Issa3}. For further information on other Hom-type algebras, one may refer to, e.g., \cite{Freg}, \cite{Goh}, \cite{Makh2}, \cite{Makh-Silv2}, \cite{Yau1}, \cite{Yau2}.
\par
The Hom-type generalization of binary algebras or ternary algebras is extended to the one of binary-ternary algebras in \cite{Gap}, \cite{Issa2}. Our present study of a Hom-type generalization of Bol algebras is included in this setting.
\par
A description of the rest of this paper is as follows.
\par
In section 2, we first recall some basics on Hom-algebras and then extend to binary-ternary Hom-algebras the notion of an $n$th-derived (binary) Hom-algebra introduced in \cite{Yau4}. Theorem 2.7 says that the category of Hom-Akivis algebras is closed under taking derived binary-ternary Hom-algebras (see Definition 2.6). Theorem 2.9, as well as Theorem 2.7, produces a sequence of Hom-Akivis. However the construction of Theorem 2.9 is not based on derived binary-ternary Hom-algebras.
\par
In section 3 we defined Hom-Bol algebras and we point out that Bol algebras are particular instances of Hom-Bol algebras. Also Hom-Bol algebras generalize multiplicative Hom-Lie triple systems in the same way as Bol algebras generalize Lie triple systems. Next we prove some construction theorems (Theorems 3.2, 3.5 and 3.6, Corollary 3.3). The category of Hom-Bol algebras is closed under self-morphisms (Theorem 3.2) and, subsequently, every Bol algebra is twisted, along any self-morphism, into a Hom-Bol algebra (Corollary 3.3). Theorem 3.5 says that the category of Hom-Bol algebras is closed under taking derived binary-ternary Hom-algebras. Theorems 3.5 and 3.6 describe some constructions of sequences of Hom-Bol algebras.
\par
In section 4, relying on a classification of real 2-dimensional Bol algebras given in \cite{Kuz-Zaid}, we classify all the algebra morphisms on all the real 2-dimensional Bol algebras and then construct (for the case of nontrivial Bol algebras) their associated Hom-Bol algebras (applying thusly Corollary 3.3).
\par
Throughout this paper we will work over a ground field of characteristic $0$.

\section{Some basics on Hom-algebras. Derived binary-ternary Hom-algebras}
The main purpose of this section is the extension to binary-ternary Hom-algebras of the notion of an $n$th derived (binary) Hom-algebra that is introduced in \cite{Yau4}. We show (Theorem 2.7) that the category of Hom-Akivis algebras is closed under taking derived binary-ternary Hom-algebras. When deriving Hom-algebras, one constructs a sequence of Hom-algebras and thus Theorem 2.7 describes a sequence of Hom-Akivis algebras. Another sequence is described in Theorem 2.9, relying on a result in \cite{Issa2}.
\par
First we begin with some basic definitions and facts that could be found in \cite{Atag}, \cite{Hart}, \cite{Issa2}, \cite{Makh2}, \cite{Makh-Silv1}, \cite{Yau1}, \cite{Yau2}, \cite{Yau5}. \\
\par
{\bf Definition 2.1.} Let $n \geq 2$ be an integer.
\par
(i) An {\it n-ary Hom-algebra} $(A, [ ,..., ], \alpha = ({\alpha}_{1},..., {\alpha}_{n-1}))$ consists of a vector space $A$, an $n$-linear map $[ ,..., ] : A^{\otimes n}  \rightarrow A$ (the $n$-ary operation) and linear maps ${\alpha}_{i} : A \rightarrow A$ (the twisting maps), $i = 1,...,n-1$.
\par
(ii) An $n$-ary Hom-algebra $(A, [ ,..., ], \alpha)$ is said to be {\it multiplicative} when the twisting maps ${\alpha}_{i}$ are all equal, ${\alpha}_{1} = ... ={\alpha}_{n-1} := \alpha$, and $\alpha ([x_{1},...,x_{n}]) = [\alpha (x_{1}),...,\alpha (x_{n})]$ for all $x_{1},...,x_{n}$ in $A$.
\par
(iii) A linear map $\theta : A \rightarrow B$ of $n$-ary Hom-algebras is called a {\it weak morphism} if $\theta ([x_{1},...,x_{n}]_{A}) = [\theta (x_{1}),...,\theta (x_{n})]_{B}$ for all $x_{1},...,x_{n}$ in $A$. The weak morphism $\theta$ is called a {\it morphism} of the $n$-ary Hom-algebras $A$ and $B$ if $\theta \circ ({\alpha}_{i})_{A} = ({\alpha}_{i})_{B} \circ \theta$ for $i = 1,...,n-1$. \\
\par
{\bf Remark.} If all $n-1$ twisting maps are the identity map $Id$ in an $n$-ary Hom-algebra Hom-algebra $(A, [ ,..., ], \alpha)$, then it reduces to an usual $n$-ary algebra $(A, [ ,..., ])$ so that the category of $n$-ary Hom-algebras contains the one of $n$-ary algebras. In this case, the weak morphism coincides with the morphism. \\
\par
For $n=2$ (resp. $n=3$), an $n$-ary Hom-algebra is called a {\it binary} (resp. {\it ternary}) Hom-algebra. In the sequel, for our purpose and convenience, we shall consider only multiplicative Hom-algebras.
\par
Hom-Lie algebras \cite{Hart} constitute the first introduced class of (binary) Hom-algebras. \\
\par
{\bf Definition 2.2.} Let $(A, [,], \alpha)$ be a multiplicative binary Hom-algebra.
\par
(i) The (binary) {\it Hom-Jacobian} of $A$ is the trilinear map $J_{\alpha} : A^{\otimes 3}  \rightarrow A$ defined as
\par
$J_{\alpha}(x,y,z) := {\circlearrowleft}_{x,y,z} [[x,y], \alpha (z)]$ \\
for all $x,y,z$ in $A$.
\par
(ii) The Hom-algebra $(A, [,], \alpha)$ is called a {\it Hom-Lie algebra} if
\par
$[x,y] = - [y,x]$,
\par
$J_{\alpha}(x,y,z) = 0$ (the {\it Hom-Jacobi identity}) \\
for all $x,y,z$ in $A$. \\
\par
We note that the notion of a Hom-Lie algebra is introduced in \cite{Hart} without the multiplicativity. Also, if $\alpha = Id$, a Hom-Lie algebra reduces to an usual Lie algebra. The $n$-ary Hom-Jacobian of an $n$-ary Hom-algebra is defined in \cite{Atag}. Other binary Hom-type algebras are introduced and discussed in \cite{Atag}, \cite{Makh1}, \cite{Makh2}, \cite{Makh-Silv1}, \cite{Yau3}, \cite{Yau4}.
\par
The classes of ternary Hom-algebras that are of interest in our setting are the ones of Hom-triple systems and Hom-Lie triple systems defined in \cite{Yau5}. \\
\par
{\bf Definition 2.3.} (i) A multiplicative {\it Hom-triple system} is a multiplicative ternary Hom-algebra $(A, [, ,], \alpha)$.
\par
(ii) A multiplicative {\it Hom-Lie triple system} is a multiplicative Hom-triple system $(A, [, ,], \alpha)$ that satisfies
\par
$[u,v,w] = -[v,u,w]$ (left skew-symmetry),
\par
${\circlearrowleft}_{u,v,w} [u,v,w] = 0$ (ternary Jacobi identity),
\par
$[\alpha (x),\alpha(y), [u,v,w]]=[[x,y,u],\alpha (v),\alpha(w)]+[\alpha (u),[x,y,v],\alpha(w)]$
\par
\hspace{3.5cm} $+ [\alpha (u), \alpha (v), [x,y,w]]$ \hfill (2.1) \\
for all $u,v,w,x,y$ in $A$. \\
\par
When $\alpha = Id$, we recover the usual notions of triple systems and Lie triple systems (\cite{Jac}, \cite{List}, \cite{Yamag}). The identity (2.1) is known as the ternary {\it Hom-Nambu identity} (see \cite{Atag} for the definition of an {\it $n$-ary Hom-Nambu algebra}). Other ternary Hom-algebras such as ternary Hom-Nambu algebras \cite{Atag}, ternary Hom-Nambu-Lie algebras \cite{Atag}, ternary Hom-Lie algebras \cite{Atag}, Hom-Jordan triple systems \cite{Yau5} are also considered. In \cite{Gap}, the Hom-triple system is defined to be a ternary Hom-algebra satsfying the left skew-symmetry and the ternary Jacobi identity (this allowed to establish a natural connection between Hom-triple systems and the Hom-version of nonassociative algebras called non-Hom-associative algebras \cite{Issa2}, or Hom-nonassociative algebras \cite{Makh2}, or either nonassociative Hom-algebras \cite{Yau2}). Hom-Lie-Yamaguti algebras introduced in \cite{Gap} constitute some generalization of multiplicative Hom-Lie triple systems in the same way as Lie-Yamaguti algebras generalize Lie triple systems. The basic object of this paper (see section 3) may also be viewed as some generalization of multiplicative Hom-Lie triple systems.
\par
Moving forward in the general theory of Hom-algebras, a study of ``binary-ternary'' Hom-algebras is initiated in \cite{Issa2} by defining the class of Hom-Akivis algebras as a Hom-analogue of the class of Akivis algebras (\cite{Ak1}, \cite{Ak2}, \cite{Hofm}) which are a typical example of binary-ternary algebras. \\
\par
{\bf Definition 2.4.} A {\it Hom-Akivis algebra} is a quadruple $(A, [,], [, ,], \alpha)$ consisting of a vector space $A$, a skew-symmetric bilinear map $[,] : A^{\otimes 2}  \rightarrow A$ (the binary operation), a trilinear map $[, ,] : A^{\otimes 3}  \rightarrow A$ (the ternary operation) and a linear self-map $\alpha$ of $A$ such that \\
\par
${\circlearrowleft}_{x,y,z} [[x,y], \alpha (z)] = {\circlearrowleft}_{x,y,z} [x,y,z] - {\circlearrowleft}_{x,y,z} [y,x,z]$  \hfill (2.2) \\
\\
for all $x,y,z$ in $A$. A Hom-Akivis algebra is said to be {\it multiplicative} if $\alpha$ is a weak morphism with respect to the binary operation ``$[, ]$'' and the ternary operation ``$[, ,]$''. \\
\par
The identity (2.2) is called the {\it Hom-Akivis identity} (if $\alpha = Id$ in (2.2), one gets the {\it Akivis identity} which defines Akivis algebras). From the definition above it clearly follows that the category of Hom-Akivis algebras contains the ones of Akivis algebras and Hom-Lie algebras. Some construction theorems for Hom-Akivis algebras are proved in \cite{Issa2}; in particular, it is shown that every Akivis algebra can be twisted along a linear self-morphism into a (multiplicative) Hom-Akivis algebra.
\par
Another class of binary-ternary Hom-algebras is the one of Hom-Lie-Yamaguti algebras \cite{Gap}.
\par
For binary Hom-algebras the notion of an $n$th derived Hom-algebra is introduced and studied in \cite{Yau4}. For completeness, we remind it in the following \\
\par
{\bf Definition 2.5.} Let $(A, \mu , \alpha)$ be a Hom-algebra and $n \geq 0$ an integer ($\mu$ is the binary operation on $A$). The Hom-algebra $A^{n}$ defined by \\
\par
$A^{n} := (A, {\mu}^{(n)} , {\alpha}^{2^n})$, \\
\\
where ${\mu}^{(n)}(x,y) := {\alpha}^{{2^n}-1}(\mu (x,y))$, $\forall x,y \in A$, is called the $n$th derived Hom-algebra of $A$. \\
\par
For simplicity of exposition, ${\mu}^{(n)}$ is written as ${\mu}^{(n)} = {\alpha}^{{2^n}-1} \circ \mu$. Then one notes that $A^{0} = (A, \mu , \alpha)$, $A^{1} = (A, {\mu}^{(1)} = \alpha \circ \mu , {\alpha}^{2})$, and $A^{n+1} = (A^{n})^1$. \\
\par
Now we extend this notion of $n$th derived (binary) Hom-algebra to the case of binary-ternary Hom-algebras in the following \\
\par
{\bf Definition 2.6.} Let ${\cal A} := (A, [,], [, ,], \alpha)$ be a binary-ternary Hom-algebra and $n \geq 0$ an integer. Define on $A$ the $n$th derived binary operation ``$[,]^{(n)}$'' and the $n$th derived ternary operation ``$[, ,]^{(n)}$'' by \\
\par
$[x,y]^{(n)} := {\alpha}^{{2^n}-1}([x,y])$, \hfill (2.3)
\par
$[x,y,z]^{(n)} := {\alpha}^{{2^{n+1}}-2}([x,y,z])$, \hfill (2.4) \\
\\
for all $x,y,z$ in $A$. Then ${\cal A}^{(n)} := (A, [,]^{(n)}, [, ,]^{(n)}, {\alpha}^{2^n})$ will be called the $n$th derived (binary-ternary) Hom-algebra of $\cal A$. \\
\par
Denote $[,]^{(n)} = {\alpha}^{{2^n}-1} \circ [,]$ and $[, ,]^{(n)} = {\alpha}^{{2^{n+1}}-2} \circ [, ,]$. Then we note that ${\cal A}^{(0)} = \cal A$, ${\cal A}^{(1)} = (A, [,]^{(1)} = \alpha \circ [,], [, ,]^{(1)} = {\alpha}^{2} \circ [, ,] , {\alpha}^{2})$ and ${\cal A}^{(n+1)} = ({\cal A}^{(n)})^{(1)}$.
\par
One observes that, from Definition 2.6, if set $[x,y,z] = 0$, $\forall x,y,z \in A$, we recover the $n$th derived (binary) Hom-algebra of Definition 2.5.
\par
The category of Hom-Akivis algebras is closed under taking derived binary-ternary Hom-algebras as it could be seen from the following result. \\
\par
{\bf Theorem 2.7.} {\it Let ${\cal A} := (A, [,], [, ,], \alpha)$ be a multiplicative Hom-Akivis algebra. Then, for each $n \geq 0$, the $n$th derived Hom-algebra ${\cal A}^{(n)} := (A, [,]^{(n)} = {\alpha}^{{2^n}-1} \circ [,], [, ,]^{(n)} = {\alpha}^{{2^{n+1}}-2} \circ [, ,], {\alpha}^{2^n})$ is a multiplicative Hom-Akivis algebra}. \\
\par
{\it Proof}. For $n=1$, ${\cal A}^{(1)} = (A, [,]^{(1)} = \alpha \circ [,], [, ,]^{(1)} = {\alpha}^{2} \circ [, ,] , {\alpha}^{2})$ is a multiplicative Hom-Akivis algebra by Theorem 4.4 in \cite{Issa2} (when $\beta = \alpha$ and $n=1$).
\par
Now suppose that, up to $n$, the ${\cal A}^{(n)}$ are multiplicative Hom-Akivis algebras. To conclude by the induction argument, we must prove that ${\cal A}^{(n+1)}$ is also a multiplicative Hom-Akivis algebra.
\par
The skew-symmetry of $[,]^{(n+1)}$ is quite obvious. In the transformations below we shall use the identity\\
\par
$J_{{\cal A}^{(n)}}(x,y,z) = {\alpha}^{2({2^n}-1)}(J_{\cal A}(x,y,z))$ \hfill (2.5) \\
\\
that holds for all $n \geq 0$ in derived (binary) Hom-algebras (see \cite{Yau4}, Lemma 2.9). We have
\begin{eqnarray}
{\circlearrowleft}_{x,y,z} [[x,y]^{(n+1)}, {\alpha}^{{2^{n+1}}} (z)]^{(n+1)} &=& J_{{\alpha}^{{2^{n+1}}}}(x,y,z) \nonumber \\
&=& {\alpha}^{2(2^{n+1}-1)}(J_{\alpha}(x,y,z)) \; \mbox{(by (2.5))} \nonumber \\
&=& {\alpha}^{2(2^{n+1}-1)}({\circlearrowleft}_{x,y,z}[x,y,z] - {\circlearrowleft}_{x,y,z}[y,x,z]) \;\; \mbox{(by (2.2))} \nonumber \\
&=& {\circlearrowleft}_{x,y,z}{\alpha}^{2^{n+2}-2}([x,y,z]) - {\circlearrowleft}_{x,y,z}{\alpha}^{2^{n+2}-2}([y,x,z]) \nonumber \\
&=& {\circlearrowleft}_{x,y,z}[x,y,z]^{(n+1)} - {\circlearrowleft}_{x,y,z}[y,x,z]^{(n+1)} \; \; \mbox{(by (2.4))} \nonumber
\end{eqnarray}
and so the Hom-Akivis identity (2.2) holds in ${\cal A}^{(n+1)}$. Thus we conclude that ${\cal A}^{(n+1)}$ is a (multiplicative) Hom-Akivis algebra. \hfill $\square$ \\
\par
As for Akivis algebras, Hom-flexibility and Hom-alternativity are defined for Hom-Akivis algebras \cite{Issa2}: a Hom-Akivis algebra ${\cal A} := (A, [,], [, ,], \alpha)$ is said to be {\it Hom-flexible} if $[x,y,x]=0$, for all $x,y$ in $A$; it is said to be {\it Hom-alternative} if $[x,y,z]=0$ whenever any two of variables $x,y,z$ are equal (i.e. $[, ,]$ is alternating). The following result constitutes the analogue of Theorem 2.7. \\
\par
{\bf Proposition 2.8.} {\it Let ${\cal A} := (A, [,], [, ,], \alpha)$ be a Hom-Akivis algebra.
\par
(i) If $\cal A$ is Hom-flexible, then the derived Hom-algebra ${\cal A}^{(n)}$ is also Hom-flexible for each $n \geq 0$.
\par
(ii) If $\cal A$ is Hom-alternative, then so is the derived Hom-algebra ${\cal A}^{(n)}$ for each $n \geq 0$.} \\
\par
{\it Proof.} We get (i) (resp. (ii)) from the definition of $[, ,]^{(n)}$ and the Hom-flexibility (resp. the Hom-alternativity) of $\cal A$. \hfill $\square$ \\
\par
Theorem 2.7 describes a sequence of Hom-Akivis algebras (the derived Hom-Akivis algebras). Starting from a given Akivis algebra and its linear self-morphism, a sequence of Hom-Akivis algebras is constructed in \cite{Issa2} (Theorem 4.8 which is a variant of Theorem 4.4). From Theorem 4.4 in \cite{Issa2} when $\alpha = \beta$, we get the following sequence of Hom-Akivis algebras but starting from a given Hom-Akivis algebra. \\
\par
{\bf Theorem 2.9.} {\it Let ${\cal A} := (A, [,], [, ,], \beta)$ be a Hom-Akivis algebra. For each integer $n \geq 0$, define on $A$ a binary operation $[,]_n$ and a ternary operation $[, ,]_n$ by 
\par
$[x,y]_{n} := {\beta}^{n}([x,y])$ (i.e. $[,]_{n} := {\beta}^{n} \circ [,]$),
\par
$[x,y,z]_{n} := {\beta}^{2n}([x,y,z])$ (i.e. $[, ,]_{n} := {\beta}^{2n} \circ [, ,]$). \\
Then ${\cal A}_{n} := (A, [,]_{n} := {\beta}^{n} \circ [,], [, ,]_{n} := {\beta}^{2n} \circ [, ,], {\beta}^{n+1})$ is a Hom-Akivis algebra}. \\
\par
{\it Proof.} The skew-symmetry of $[,]_n$ is obvious. Next, we have
\begin{eqnarray}
{\circlearrowleft}_{x,y,z} [[x,y]_{n}, {\beta}^{n+1} (z)]_{n} &=& {\circlearrowleft}_{x,y,z} {\beta}^{n}([{\beta}^{n}([x,y]), {\beta}^{n+1}(z)]) \nonumber \\
&=& {\beta}^{2n}({\circlearrowleft}_{x,y,z}[[x,y], {\beta} (z)]) \nonumber \\
&=& {\beta}^{2n}({\circlearrowleft}_{x,y,z}[x,y,z] - {\circlearrowleft}_{x,y,z}[y,x,z]) \; \mbox{(by (2.2) in} \; \cal A) \nonumber \\
&=& {\circlearrowleft}_{x,y,z}{\beta}^{2n}([x,y,z]) - {\circlearrowleft}_{x,y,z}{\beta}^{2n}([y,x,z]) \nonumber \\
&=& {\circlearrowleft}_{x,y,z}[x,y,z]_{n} - {\circlearrowleft}_{x,y,z}[y,x,z]_{n} \nonumber
\end{eqnarray}
which means that (2.2) holds in ${\cal A}_{n}$ and so ${\cal A}_{n}$ is a Hom-Akivis algebra. \hfill $\square$ \\
\par
In Theorem 2.9 observe that ${\cal A}_{0} = \cal A$. However, in strong contrast with derived Hom-Akivis algebras, ${\cal A}_{n+1}$ is not constructed from ${\cal A}_{n}$.
\section{Hom-Bol algebras. Definition and construction theorems}
In this section we define a Hom-Bol algebra. It turns out that Hom-Bol algebras constitute a twisted generalization of Bol algebras. We prove some construction theorems for Hom-Bol algebras (Theorems 3.2, 3.5, Corollary 3.3). Theorem 3.2 shows that the category of Hom-Bol algebras is closed under self-morphisms while Corollary 3.3 says that every Bol algebra can be twisted, along any endomorphism, into a Hom-Bol algebra. Theorem 3.5 points out that Hom-Bol algebras are also closed under taking derived binary-ternary Hom-algebras.
\par
We begin with the definition of the basic object of this paper. \\
\par
{\bf Definition 3.1.} A {\it Hom-Bol algebra} is a quadruple $(A,*, \{ , , \}, \alpha)$ in which $A$ is a vector space, ``$*$'' a binary operation and ``$\{ , , \}$'' a ternary operation on $A$, and $\alpha : A \rightarrow A$ a linear map such that\\
\par
(HB1) $\alpha ( x * y) = \alpha (x) * \alpha (y)$,
\par
(HB2) $\alpha (\{x,y,z\}) = \{\alpha (x), \alpha (y), \alpha (z) \}$,
\par
(HB3) $x*y = - y*x$,
\par
(HB4) $\{x,y,z\} = - \{y,x,z\}$,
\par
(HB5) ${\circlearrowleft}_{x,y,z} \{x,y,z\} = 0$,
\par
(HB6) $\{\alpha (x), \alpha (y), u*v \} = \{x,y,u\}* {\alpha}^{2}(v) + {\alpha}^{2}(u)*\{x,y,v\} + \{\alpha (u), \alpha (v), x*y \}$
\par
\hspace{4.5cm}$- \alpha (u)\alpha (v) * \alpha (x)\alpha (y)$,
\par
(HB7) $\{ {\alpha}^{2}(x), {\alpha}^{2}(y), \{ u,v,w \} \} = \{ \{x,y,u\}, {\alpha}^{2}(v), {\alpha}^{2}(w) \} + \{ {\alpha}^{2}(u), \{x,y,v\}, {\alpha}^{2}(w) \}$
\par
\hspace{5.5cm}$+ \{ {\alpha}^{2}(u), {\alpha}^{2}(v), \{x,y,w\} \}$\\
\\
for all $u,v,w,x,y,z$ in $A$. \\
\par
The multiplicativity of $(A,*, \{ , , \}, \alpha)$ (see (HB1) and (HB2)) is built into our definition for convenience. \\
\par
{\bf Remark.} (i) If $\alpha = Id$, then the Hom-Bol algebra $(A,*, \{ , , \}, \alpha)$ reduces to a Bol algebra $(A,*, \{ , , \})$ (see (B1)-(B5)). So Bol algebras may be seen as Hom-Bol algebras with the identity map $Id$ as the twisting map.
\par
(ii) If $x*y = 0$, for all $x,y \in A$, then $(A,*, \{ , , \}, \alpha)$ becomes a (multiplicative) Hom-Lie triple system $(A, \{ , , \}, {\alpha}^{2})$ (see Definition 2.3). The identity (HB7) is in fact the ternary Hom-Nambu identity (2.1) but with ${\alpha}^{2}$ as the twisting map.\\
\par
The following result shows that the category of Hom-Bol algebras is closed under self-morphisms. \\
\par
{\bf Theorem 3.2.} {\it Let ${\cal A}_{\alpha} := (A,*, \{ , , \}, \alpha)$ be a Hom-Bol algebra and $\beta$ and endomorphism of $(A,*, \{ , , \})$ such that $\beta \alpha = \alpha \beta$ (i.e. $\beta$ is a self-morphism of ${\cal A}_{\alpha}$). Let ${\beta}^{0} = Id$ and, for any $n \geq 1$, ${\beta}^{n} := \beta {\beta}^{n-1} \; (= \beta \circ {\beta}^{n-1})$. Define on $A$ the operations
\par
$ x *_{\beta} y := {\beta}^{n} (x * y)$ (i.e. $*_{\beta} := {\beta}^{n} \circ *$),
\par
$\{x,y,z \}_{\beta} := {\beta}^{2n}(\{x,y,z \})$ (i.e. $\{ ,, \}_{\beta} := {\beta}^{2n} \circ \{ ,, \}$)\\
for all $x,y,z$ in $A$. Then ${\cal A}_{{\beta}^{n}} := (A, *_{\beta}, \{ ,, \}_{\beta}, {\beta}^{n} \alpha)$ is a Hom-Bol algebra with $n \geq 1$.}\\
\par
{\it Proof.} We observe first that the condition $\beta \alpha = \alpha \beta$ implies ${\beta}^{n} \alpha = \alpha {\beta}^{n}$. This last equality and the definitions of the operations ``$*_{\beta}$'' and ``$\{ ,, \}_{\beta}$'' lead to the validity of (HB1) and (HB2) for ${\cal A}_{{\beta}^{n}}$. Obviously the skew-symmetries (HB3) and (HB4) for ${\cal A}_{{\beta}^{n}}$ follow from the skew-symmetry of ``$*$'' and ``$\{ ,, \}$'' respectively. Next, we have
\par
${\circlearrowleft}_{x,y,z} \{x,y,z \}_{\beta} = {\circlearrowleft}_{x,y,z} {\beta}^{2n} (\{x,y,z \}) = {\beta}^{2n} ({\circlearrowleft}_{x,y,z} (\{ x, y, z \}) = {\beta}^{2n}(0)$ (by (HB5) for ${\cal A}_{\alpha}$)
\par
$=0$\\
so we get (HB5) for ${\cal A}_{{\beta}^{n}}$.
\par
Consider now $\{({\beta}^{n} \alpha)(x), ({\beta}^{n} \alpha)(y), u \; *_{\beta} \; v \}_{\beta}$ in ${\cal A}_{{\beta}^{n}}$. We have (using the condition $\beta \alpha = \alpha \beta$)
\begin{eqnarray}
& & \{({\beta}^{n} \alpha)(x), ({\beta}^{n} \alpha)(y), u \; *_{\beta} \; v \}_{\beta}={\beta}^{3n}(\{ \alpha (x), \alpha (y), u*v \}) \nonumber \\
&=& {\beta}^{3n}( \{ x, y, u \} * {\alpha}^{2}(v) + {\alpha}^{2}(u) * \{ x, y, v \} + \{\alpha (u), \alpha (v), x*y \} \nonumber \\
&-& \alpha (u)\alpha (v) * \alpha (x)\alpha (y)) \; \mbox{(by (HB6) for} \; (A,*, \{ , , \}, \alpha)) \nonumber \\
&=& {\beta}^{n}({\beta}^{2n}(\{ x, y, u \})*({\beta}^{2n}{\alpha}^{2})(v) + {\beta}^{n}(({\beta}^{2n}{\alpha}^{2})(u)*{\beta}^{2n}(\{ x, y, v \})) \nonumber \\
&+& {\beta}^{2n}( \{ ({\beta}^{n} \alpha)(u), ({\beta}^{n} \alpha)(v), {\beta}^{n}(x*y) \}) - {\beta}^{n} ({\beta}^{2n}(\alpha (u)\alpha (v)) * {\beta}^{2n}(\alpha (x)\alpha (y))) \nonumber \\
&=&\{ x,y,u \}_{\beta} *_{\beta} ({\beta}^{n} \alpha)^{2}(v) + ({\beta}^{n} \alpha)^{2}(u) *_{\beta} \{ x,y,v \}_{\beta} \nonumber \\
&+& \{ ({\beta}^{n} \alpha)(u), ({\beta}^{n} \alpha)(v), x *_{\beta} y \}_{\beta} - (({\beta}^{n} \alpha)(u) *_{\beta} ({\beta}^{n} \alpha)(v)) *_{\beta} (({\beta}^{n} \alpha)(x) *_{\beta} ({\beta}^{n} \alpha)(y)) \nonumber
\end{eqnarray}
and thus we get (HB6) for ${\cal A}_{{\beta}^{n}}$. Finally, using repeatedly the condition $\beta \alpha = \alpha \beta$ and (HB7) for ${\cal A}_{\alpha}$, we compute:
\begin{eqnarray}
& & \{ {({\beta}^{n} \alpha)}^{2}(x), {({\beta}^{n} \alpha)}^{2}(y), \{ u, v, w \}_{\beta} \}_{\beta} \nonumber \\
&=& \{({\beta}^{2n} {\alpha}^{2})(x), ({\beta}^{2n} {\alpha}^{2})(y), \{ u, v, w\}_{\beta} \}_{\beta} \nonumber \\
&=&  {\beta}^{2n}(\{ ({\beta}^{2n} {\alpha}^{2})(x), ({\beta}^{2n} {\alpha}^{2})(y), {\beta}^{2n}(\{ u, v, w \})\}) \nonumber \\
&=& {\beta}^{4n}(\{ {\alpha}^{2}(x), {\alpha}^{2}(y), \{ u, v, w \} \}) \nonumber \\
&=& {\beta}^{4n}(\{{\alpha}^{2}(u), {\alpha}^{2}(v), \{x,y,w \} \}) \nonumber \\
&+& {\beta}^{4n}(\{ \{x,y,u \}, {\alpha}^{2}(v), {\alpha}^{2}(w) \}) \nonumber \\
&+& {\beta}^{4n}(\{{\alpha}^{2}(u), \{x,y,v\}, {\alpha}^{2}(w) \} \nonumber \\
&=&  {\beta}^{2n}(\{({\beta}^{2n} {\alpha}^{2})(u), ({\beta}^{2n} {\alpha}^{2})(v), {\beta}^{2n}(\{x,y,w \})\}) \nonumber \\
&+& {\beta}^{2n}(\{ {\beta}^{2n}(\{x,y,u\}), ({\beta}^{2n} {\alpha}^{2})(v), ({\beta}^{2n} {\alpha}^{2})(w)\}) \nonumber \\
&+& {\beta}^{2n}(\{({\beta}^{2n} {\alpha}^{2})(u),{\beta}^{2n}(\{x,y,v\}), ({\beta}^{2n} {\alpha}^{2})(w)\}) \nonumber \\
&=& \{ {({\beta}^{n} \alpha)}^{2}(u), {({\beta}^{n} \alpha)}^{2}(v), \{x,y,w \}_{\beta} \}_{\beta} \nonumber \\
&+& \{ \{x,y,u \}_{\beta}, {({\beta}^{n} \alpha)}^{2}(v), {({\beta}^{n} \alpha)}^{2}(w) \}_{\beta} \nonumber \\
&+& \{ {({\beta}^{n} \alpha)}^{2}(u),\{x,y,v \}_{\beta}, {({\beta}^{n} \alpha)}^{2}(w) \}_{\beta}. \nonumber
\end{eqnarray}
Thus (HB7) holds for ${\cal A}_{{\beta}^{n}}$.
\par
Therefore, we proved the validity for ${\cal A}_{{\beta}^{n}}$ of the set of identities of type (HB1)-(HB7) and so we get that ${\cal A}_{{\beta}^{n}}$ is a Hom-Bol algebra. This finishes the proof. \hfill $\square$\\
\par
From Theorem 3.2 we have the following\\
\par
{\bf Corollary 3.3.} {\it Let ${\cal A} := (A,*, [, ,])$ be a Bol algebra and $\beta$ an endomorphism of $\cal A$. If define on $A$ a binary operation ``$ \widetilde{*}$'' and a ternary operation ``$\{ , , \}$'' by
\par
$x \widetilde{*} y := \beta (x*y)$,
\par
$\{x,y,z\} := {\beta}^{2} ([x,y,z])$, \\
then $(A, \widetilde{*} , \{ , , \}, \beta)$ is a Hom-Bol algebra}. \\
\par
{\it Proof.} The proof follows from the one of Theorem 3.2 for $\alpha = Id$ and $n=1$. \hfill $\square$ \\
\par
Observe that Corollary 3.3 gives a method for constructing Hom-Bol algebras from Bol algebras. This is an extension of a result due to D. Yau \cite{Yau2} giving a general construction method of Hom-algebras from their corresponding untwisted version (see also \cite{Gap}, \cite{Issa2} for a use of such an extension). From Corollary 3.3, we get the following \\
\par
{\bf Proposition 3.4.} {\it Let $(A,*)$ be a Malcev algebra and $\beta$ any endomorphism of $(A,*)$. Define on $A$ the operations}
\par
$x \tilde{*} y := \beta (x * y)$,
\par
$\{ x,y,z \} := (1/3){\beta}^{2}(2(x*y)*z - (y*z)*x - (z*x)*y)$. \\
{\it Then $(A, \tilde{*}, \{ , , \}, \beta)$ is a Hom-Bol algebra}. \\
\par
{\bf Proof.} If consider on $A$ the ternary operation $[x,y,z] :=$ \\ $(1/3)(2(x*y)*z - (y*z)*x - (z*x)*y)$, $\forall x,y,z \in A$, then $(A, * , [ , , ])$ is a Bol algebra \cite{Mikh}. Moreover, since $\beta$ is an endomorphism of $(A,*)$, we have $\beta ([x,y,z]) = [\beta(x),\beta(y),\beta(z)]$ so that $\beta$ is also an endomorphism of $(A, * , [ , , ])$. Then Corollary 3.3 implies that $(A, \tilde{*}, \{ , , \}, \beta)$ is a Hom-Bol algebra. \hfill $\square$ \\
\par
The following construction result for Hom-Bol algebras is the analogue of Theorem 2.7. It shows that Hom-Bol algebras are closed under taking derived Hom-algebras. \\
\par
{\bf Theorem 3.5.} {\it Let ${\cal A} := (A,*, \{ , , \}, \alpha)$ be a Hom-Bol algebra. Then, for each $n \geq 0$, the $n$th derived Hom-algebra ${{\cal A}^{(n)}} := (A, *^{(n)} = {\alpha}^{2^{n}-1} \circ *, \{, ,\}^{(n)} = {\alpha}^{{2^{n+1}}-2} \circ \{, ,\}, {\alpha}^{2^n})$ is a Hom-Bol algebra}. \\
\par
{\it Proof.} The identities (HB1)-(HB5) for ${\cal A}^{(n)}$ are obvious. The checking of (HB6) for ${\cal A}^{(n)}$ is as follows.
\begin{eqnarray}
\{ {\alpha}^{2^n}(x), {\alpha}^{2^n}(y), u *^{(n)} v\}^{(n)} &=& {\alpha}^{{2^{n+1}}-2}(\{{\alpha}^{2^n}(x), {\alpha}^{2^n}(y), {\alpha}^{2^{n-1}}(u * v) \}) \nonumber \\
&=& {\alpha}^{{2^{n+1}}-2}(\{ {\alpha}^{2^{n-1}} \alpha (x), {\alpha}^{2^{n-1}} \alpha (y), {\alpha}^{2^{n-1}}(u * v) \}) \nonumber \\
&=& {\alpha}^{{2^{n+1}}-2} {\alpha}^{2^{n-1}} (\{ \alpha (x), \alpha (y), u * v \}) \nonumber \\
&=& {\alpha}^{{2^{n+1}}-2} {\alpha}^{2^{n-1}} (\{ x,y,u \} * {\alpha}^{2}(v) + {\alpha}^{2}(u) * \{ x,y,v \}\nonumber \\
&+& \{ \alpha (u), \alpha (v), x * y \} - \alpha (u) \alpha (v) * \alpha (x) \alpha (y)) \nonumber \\
& & \mbox{(by (HB6) for} \; \cal A) \nonumber \\
&=& {\alpha}^{{2^{n+1}}-2}( \{ x,y,u \} *^{(n)} {\alpha}^{2}(v)) \nonumber \\
&+& {\alpha}^{{2^{n+1}}-2} ({\alpha}^{2}(u) *^{(n)} \{ x,y,v \}) \nonumber \\
&+& {\alpha}^{{2^{n+1}}-2} (\{ {\alpha}^{2^n}(u), {\alpha}^{2^n}(v), x *^{(n)} y \}) \nonumber \\
&-& {\alpha}^{{2^{n+1}}-2} (\alpha (u) \alpha (v) *^{(n)} \alpha (x) \alpha (y)) \nonumber \\
&=& \{ x,y,u \}^{n}  *^{(n)} {\alpha}^{2^{n+1}} (v) + {\alpha}^{2^{n+1}} (u) *^{(n)} \{ x,y,v \}^{n} \nonumber \\
&+& \{ {\alpha}^{2^n}(u), {\alpha}^{2^n}(v), x *^{(n)} y \}^{(n)} \nonumber \\
&-& ({\alpha}^{2^{n-1}})^{2}(\alpha (u) \alpha (v) *^{(n)} \alpha (x) \alpha (y)) \nonumber \\
&=& \{ x,y,u \}^{n}  *^{(n)} ({\alpha}^{2^n})^{2} (v) + ({\alpha}^{2^n})^{2} (u) *^{(n)} \{ x,y,v \}^{n} \nonumber \\
& + & \{ {\alpha}^{2^n}(u), {\alpha}^{2^n}(v), x *^{(n)} y \}^{(n)} \nonumber \\
& - & ({\alpha}^{2^n}(u) *^{(n)} {\alpha}^{2^n}(v)) *^{(n)} ({\alpha}^{2^n}(x) *^{(n)} {\alpha}^{2^n}(y)) \nonumber
\end{eqnarray}
and thus (HB6) holds for ${\cal A}^{(n)}$. Finally, we compute
\begin{eqnarray}
& & \{ ({\alpha}^{2^n})^{2}(x), ({\alpha}^{2^n})^{2}(y), \{ u,v,w \}^{(n)} \}^{(n)} \nonumber \\
&=& {\alpha}^{{2^{n+1}}-2} (\{ ({\alpha}^{2^n})^{2}(x), ({\alpha}^{2^n})^{2}(y), {\alpha}^{{2^{n+1}}-2}(\{ u,v,w \})\}) \nonumber \\
&=& {\alpha}^{{2^{n+1}}-2} (\{ {\alpha}^{2^{n+1}} (x), {\alpha}^{2^{n+1}} (y), {\alpha}^{{2^{n+1}}-2} (\{u,v,w \})\}) \nonumber \\
&=& ({\alpha}^{{2^{n+1}}-2})^{2} ( \{ {\alpha}^{2}(x), {\alpha}^{2}(y), \{u,v,w \} \} ) \nonumber \\
&=& ({\alpha}^{{2^{n+1}}-2})^{2} ( \{\{ x,y,u \}, {\alpha}^{2}(v), {\alpha}^{2}(w) \}) + ({\alpha}^{{2^{n+1}}-2})^{2} (\{ {\alpha}^{2}(u), \{ x,y,v \}, {\alpha}^{2}(w) \}) \nonumber \\
&+& ({\alpha}^{{2^{n+1}}-2})^{2} (\{ {\alpha}^{2}(u), {\alpha}^{2}(v), \{ x,y,w \} \}) \; \mbox{(by (HB7) for} \; \cal A) \nonumber \\
&=& {\alpha}^{{2^{n+1}}-2} (\{ {\alpha}^{{2^{n+1}}-2} (\{ x,y,u \}), {\alpha}^{{2^{n+1}}} (v), {\alpha}^{{2^{n+1}}} (w) \}) \nonumber \\
&+& {\alpha}^{{2^{n+1}}-2} (\{{\alpha}^{{2^{n+1}}} (u), {\alpha}^{{2^{n+1}}-2} (\{ x,y,v \}), {\alpha}^{{2^{n+1}}} (w) \}) \nonumber \\
&+& {\alpha}^{{2^{n+1}}-2} (\{{\alpha}^{{2^{n+1}}} (u), {\alpha}^{{2^{n+1}}} (v), {\alpha}^{{2^{n+1}}-2} (\{ x,y,v \})\}) \nonumber \\
&=& \{\{ x,y,u \}^{(n)}, ({\alpha}^{2^n})^{2}(v), ({\alpha}^{2^n})^{2}(w) \}^{(n)} + \{ ({\alpha}^{2^n})^{2}(u), \{ x,y,v \}^{(n)}, ({\alpha}^{2^n})^{2}(w) \}^{(n)} \nonumber \\
&+& \{ ({\alpha}^{2^n})^{2}(u), ({\alpha}^{2^n})^{2}(v), \{ x,y,w \}^{(n)} \}^{(n)} \nonumber
\end{eqnarray}
and so we see that (HB7) holds for ${\cal A}^{(n)}$. Thus we get that ${\cal A}^{(n)}$ verifies the system (HB1)-(HB7), which means that ${\cal A}^{(n)}$ is a Hom-Bol algebra as claimed. \hfill $\square$ \\
\par
The construction described in Theorem 2.9 is reported to Hom-Bol algebras as follows. \\
\par
{\bf Theorem 3.6.} {\it Let ${\cal A} := (A,*, \{ , , \}, \beta)$ be a Hom-Bol algebra. For each $n \geq 0$, the Hom-algebra ${\cal A}_{n} := (A, *_{n}={\beta}^{n} \circ * , \{ , , \}_{n}= {\beta}^{2n} \circ \{ , , \}, {\beta}^{n+1})$ is a Hom-Bol algebra.} \\
\par
{\it Proof.} We get the proof from the one of Theorem 3.2 if set $\alpha = \beta$ in Theorem 3.2. \hfill $\square$ \\
\section{Examples of 2-dimensional Hom-Bol algebras}
The purpose of this section is to classify all the algebra morphisms on all the 2-dimensional Bol algebras, using the classification of real 2-dimensional Bol algebras given in \cite{Kuz-Zaid} (another classification, infered from the one of 2-dimensional hyporeductive triple algebras, is independently given in \cite{Issa1}, but not as distinct isomorphic classes). Then, from Corollary 3.3 we obtain their corresponding 2-dimensional Hom-Bol algebras. In this section we will work over the ground field of real numbers.
\par
We consider only nontrivial cases, i.e. we will omit mentioning the Lie triple systems (as particular cases of Bol algebras) and the 0-map as twisting map (since it is always an algebra morphism and it gives rise to the zero Hom-algebra).
\par
If $(A, *, [ , , ])$ is a 2-dimensional Bol algebra with basis $\{ e_{1}, e_{2} \}$, then a linear map $\beta : A \rightarrow A$ defined by $\beta (e_{1}) = a_{1} e_{1} + a_{2} e_{2}$, $\beta (e_{2}) = b_{1} e_{1} + b_{2} e_{2}$ is a self-morphism of $(A, *, [ , , ])$ if and only if\\
\par
$\beta (e_{1} * e_{2}) = \beta (e_{1}) * \beta (e_{2})$, $\beta ([e_{1}, e_{1}, e_{i}]) = [\beta (e_{1}), \beta (e_{2}), \beta (e_{i})]$, $i=1,2$. \hfill (4.1) \\
\par
It is proved (\cite{Kuz-Zaid}, Theorem 2) that any nontrivial 2-dimensional (left) Bol algebra over the field of real numbers is isomorphic to only one of the algebras described below:\\
\par
(A1) $e_{1} * e_{2} = -e_{2}$, $[e_{1},e_{2},e_{1}] = e_{1}$, $ [e_{1},e_{2},e_{2}] = -e_{2}$;
\par
(A2) $e_{1} * e_{2} = -e_{2}$, $[e_{1},e_{2},e_{1}] = \lambda e_{2}$, $ [e_{1},e_{2},e_{2}] =0$;
\par
(A3) $e_{1} * e_{2} = -e_{2}$, $[e_{1},e_{2},e_{1}] = \lambda e_{2}$, $[e_{1},e_{2},e_{2}] = \pm e_{1}$. \\
\par
We recall that we consider here left Bol algebras, instead of the right-sided case studied in \cite{Kuz-Zaid}. Also, in \cite{Kuz-Zaid} the types (A2) and (A3) are gathered in a single one type.
\par
Now suppose that, with respect to the basis $\{ e_{1}, e_{2} \}$, the linear map $\beta$ is given as $\beta (e_{1}) = a_{1} e_{1} + a_{2} e_{2}$, $\beta (e_{2}) = b_{1} e_{1} + b_{2} e_{2}$. We shall discuss the conditions (4.1) for each of the type (A1), (A2) and (A3). \\
\par
$\star$ Case of type (A1):
\par
For this type, the conditions (4.1) lead to the following simultaneous constraints on the coefficients $a_{i}$, $b_{j}$ of $\beta$: \\
\begin{displaymath}
\left\{ \begin{array}{ll}
b_{1} = 0 \\
b_{2} (a_{1} -1) = 0,
\end{array} \right.
\end{displaymath}
\begin{displaymath}
\left\{ \begin{array}{ll}
a_{1}(1-a_{1} b_{2})= 0 \\
a_{2}(1 + a_{1}b_{2}) = 0,
\end{array} \right.
\end{displaymath}
$$b_{2}(1-a_{1} b_{2})= 0.$$
\par
From these constraints, we see that the identity map $\beta = Id$ is the only nonzero self-morphism of (A1). So, by Corollary 3.3, any Bol algebra of type (A1), regarded as a Hom-Bol algebra with $Id$ as the twisting map, is the only one nonzero Hom-Bol algebra correspondig to (A1). \\
\par
$\star$ Case of type (A2):
\par
The conditions (4.1) then imply that
\begin{displaymath}
\left\{ \begin{array}{ll}
b_{1} = 0 \\
b_{2} (a_{1} -1) = 0,
\end{array} \right.
\end{displaymath}
$$\lambda b_{2}(1 - {a_{1}}^{2}) =0.$$ \\
\par
- If $b_{2} = 0$, from the equations above, we see that the nonzero morphism $\beta$ is defined as \\
\par
$\beta (e_{1}) = a_{1} e_{1} + a_{2} e_{2}$, (with $(a_{1},a_{2}) \neq (0,0)$), $\beta (e_{2}) = 0$ \hfill (4.2) \\
\\
and, applying Corollary 3.3, we get that (A2) is twisted into the zero Hom-Bol algebra.
\par
- If $b_{2} \neq 0$, then the equation $b_{2} (a_{1} -1) = 0$ (see above) implies $a_{1}=1$ and so the nonzero morphism $\beta$ is defined as \\
\par
$\beta (e_{1}) = e_{1} + a_{2} e_{2}$, $\beta (e_{2}) = b_{2} e_{2}$ \; ($b_{2} \neq 0$). \hfill (4.3) \\
\par
The application of Corollary 3.3 then gives rise to the Hom-Bol algebra $(A, \widetilde{*} , \{ , , \}, \beta)$ defined by \\
\par
$e_{1} \widetilde{*} e_{2} = - b_{2} e_{2}$, $\{ e_{1}, e_{2},  e_{1} \} = \lambda b_{2} e_{2}$, $\{ e_{1}, e_{2},  e_{2} \} = 0$ \hfill (4.4) \\
\\
with $b_{2} \neq 0$. \\
\par
$\star$ Case of type (A3):
\par
In this case, the conditions (4.1) are expressed as the following simultaneous equations:
\begin{displaymath}
\left\{ \begin{array}{ll}
b_{1} = 0 \\
b_{2} (a_{1} -1) = 0,
\end{array} \right.
\end{displaymath}
\begin{displaymath}
\left\{ \begin{array}{ll}
a_{1} a_{2} b_{2} = 0 \\
\lambda b_{2}(1 - {a_{1}}^{2}) =0,
\end{array} \right.
\end{displaymath}
\begin{displaymath}
\left\{ \begin{array}{ll}
a_{1}(1 - {b_{2}}^{2}) =0\\
a_{2}= 0.
\end{array} \right.
\end{displaymath}
So $\beta$ must be found in the form\\
\par
$\beta (e_{1}) = a_{1}e_{1}$, $\beta (e_{2}) = b_{2} e_{2}$. \\
\\
-If $b_{2} = 0$, then the equation $a_{1}(1 - {b_{2}}^{2}) =0$ implies $a_{1}= 0$ and then $\beta$ is the 0-map, which is the case that is already precluded. \\
\\
- Let $b_{2} \neq 0$. Then the equation $b_{2} (a_{1} -1) = 0$ implies $a_{1}= 1$ and so $\beta$ is defined as \\
\par
$\beta (e_{1}) = e_{1}$, $\beta (e_{2}) = b_{2} e_{2}$ \; ($b_{2} \neq 0$). \hfill (4.5) \\
\par
Therefore Corollary 3.3 says that (A3) is twisted into the Hom-Bol algebra \\ $(A, \widetilde{*} , \{ , , \}, \beta)$ defined by \\
\par
$e_{1} \widetilde{*} e_{2} = - b_{2} e_{2}$, $\{ e_{1}, e_{2},  e_{1} \} = \lambda b_{2} e_{2}$, $\{ e_{1}, e_{2},  e_{2} \} = \mp e_{1}$ \hfill (4.6) \\
\\
($b_{2} \neq 0$).
\par
Summarizing the considerations above, we conclude that:
\par
(i) the algebra (A1) (regarded as a Hom-Bol algebra) is the only (up to isomorphism) nonzero Hom-Bol algebra corresponding to a Bol algebra of type (A1), the identity map $Id$ being the only nonzero self-morphism of (A1);
\par
(ii) a Bol algebra of type (A2) is twisted into either a zero Hom-Bol algebra, (by a morphism (4.2)) or a Hom-Bol algebra of type (4.4) (by a morphism (4.3));
\par
(iii) a Bol algebra of type (A3) is twisted (precluding the 0-map) into a Hom-Bol algebra of type (4.6) (by a morphism (4.5)).
\par
From the other hand, applying Theorem 3.5, we see that the Hom-Bol algebras (A1), (4.4), and (4.6) yield the following sequences of $n$th derived Hom-algebras (denoting, for simplicity, $b_{2} := b$ and $a_{2} := a$):
\par
{\textbullet} For each $n \geq 0$, the nonzero $n$th derived Hom-Bol algebra of (A1) (regarded as a Hom-Bol algebra) is (A1) itself (recall that the identity map $Id$ is the only nonzero self-morphism of (A1)).
\par
{\textbullet} For each $n \geq 0$, the $n$th derived Hom-Bol algebra of the nonzero Hom-Bol algebra (4.4) is given by\\
\par
$e_{1} {\widetilde{*}}^{(n)} e_{2} = - {b}^{2^{n}-1} e_{2}$, $\{ e_{1}, e_{2},  e_{1} \}^{(n)} = \lambda {b}^{2^{n+1}-1} e_{2}$, $\{ e_{1}, e_{2},  e_{2} \}^{(n)} = 0$ \\
\\
with the twisting map ${\beta}^{2^n}$ defined by\\
\par
${\beta}^{2^n} (e_{1}) = e_{1} + a(1+b+ ... +{b}^{2^{n}-1})e_{2}$, ${\beta}^{2^n} (e_{2}) = {b}^{2^n} e_{2}$.\\
\par
{\textbullet} For each $n \geq 0$, the $n$th derived Hom-Bol algebra of the nonzero Hom-Bol algebra (4.6) is given by\\
\par
$e_{1} {\widetilde{*}}^{(n)} e_{2} = - {b}^{2^{n}-1} e_{2}$, $\{ e_{1}, e_{2},  e_{1} \}^{(n)} = \lambda {b}^{2^{n+1}-1} e_{2}$, $\{ e_{1}, e_{2},  e_{2} \}^{(n)} = \mp e_{1}$ \\
\\
with the twisting map ${\beta}^{2^n}$ defined by\\
\par
${\beta}^{2^n} (e_{1}) = e_{1}$, ${\beta}^{2^n} (e_{2}) = {b}^{2^n} e_{2}$.

Sylvain Attan \\ Institut de Math\'ematiques et de Sciences Physiques, \\Universit\'e d'Abomey-Calavi \\
01 BP 613-Oganla, Porto-Novo, B\'enin \\ sylvain.attan@imsp-uac.org \\
\\
A. Nourou Issa \\ D\'epartement de Math\'ematiques, Universit\'e d'Abomey-Calavi, \\ 01 BP 4521, Cotonou 01, B\'enin \\ woraniss@yahoo.fr

\end{document}